\documentclass[11pt]{article}
\usepackage{latexsym}

\newcommand{\qed}{\ \hfill\mbox{$\Box$}\vspace{\baselineskip}}

\newtheorem{theorem}{Theorem}[section]
\newtheorem{lemma}{Lemma}
\newtheorem{proposition}[theorem]{Proposition}

\newenvironment{proof}{\noindent {\bf Proof:}}{{\qed}}

\newcommand{\iv}{\ensuremath{{\mathcal I}}}

\newcommand{\rank}{\mbox{rank}}
\newcommand{\supp}{\mbox{supp}}

\begin{document}

\title{Signs in the $cd$-index of Eulerian~partially~ordered~sets}

\author{Margaret M. Bayer\\
        Department of Mathematics\\
	University of Kansas\\
	Lawrence KS 66045-2142\\
        bayer@math.ukans.edu}

\date{}

\maketitle

\begin{abstract}
A graded partially ordered set is Eulerian if every interval has the same
number of elements of even rank and of odd rank.
Face lattices of convex polytopes are Eulerian.
For Eulerian partially ordered sets, the flag vector can be encoded efficiently
in the $cd$-index.
The $cd$-index of a polytope has all positive entries.
An important open problem is to give the broadest natural class of Eulerian
posets having nonnegative $cd$-index.
This paper completely determines which entries of the $cd$-index are nonnegative
for all Eulerian posets.
It also shows that there are no other lower or upper bounds on $cd$-coefficients
(except for the coefficient of $c^n$).
\end{abstract}

\maketitle

\section{Introduction}
In the past thirty years or more, there has been much interest in combinatorial
questions about polytopes and other geometric complexes and partial orders.
Of central importance is the flag vector of a partially ordered set (poset)
and various combinatorial parameters derived from it.
One of these parameters is the $cd$-index, defined for Eulerian posets, a class
that contains face lattices of polytopes.
The $cd$-index was discovered by Fine and introduced in the literature by 
Bayer and Klapper (\cite{Bayer-Klapper}).
It has captured the imagination, both for what is known and for what is not 
known about it.
It embodies in an elegant way the linear relations of flag vectors of Eulerian 
posets (the generalized Dehn-Sommerville relations of Bayer and Billera 
\cite{Bayer-Billera}); the number of coefficients in the $cd$-index is a
Fibonacci number.
It is known to be nonnegative for polytopes (see Stanley \cite{Stanley-flag}),
but it is not known what it counts, except in special cases (see Purtill
\cite{Purtill}).
Among polytopes, the $cd$-index is minimized by the simplices (see Billera and
Ehrenborg \cite{Billera-Ehrenborg}).
Novik (\cite{Novik}) gives lower bounds for $cd$-coefficients of odd-dimensional
simplicial manifolds (or, more generally, Eulerian Buchsbaum complexes).

Stanley (\cite{Stanley-flag}) proved the nonnegativity of the $cd$-index for
``S-shellable'' regular CW-spheres (including polytopes).
In \cite{Stanley-posprob} he proposes the following as the main open problem 
concerning the $cd$-index:  Is the $cd$-index nonnegative for all
Gorenstein$^*$ posets?
(These are the Cohen-Macaulay Eulerian posets.)
In fact, some parts of the $cd$-index are nonnegative for {\em all}
Eulerian posets.  
In this paper we determine which $cd$-words have nonnegative coefficients for
all Eulerian posets.
For all other $cd$-words, we show how to construct Eulerian posets with 
arbitrarily large negative coefficients.
The proofs grow out of the ideas of \cite{Bayer-Hetyei1}, which studies the
cone of flag vectors of Eulerian posets.

\section{Definitions}
An {\em Eulerian poset} is a graded partially ordered set $P$ in which every
interval has the same number of elements of even and of odd rank.
For $P$ an Eulerian poset, the dual poset, obtained by reversing the order
relation, is also Eulerian.
The $cd$-index of an Eulerian poset is an invariant based on the numbers of 
chains in the poset.
For $P$ an Eulerian poset of rank $n+1$ and $S\subseteq [1,n]$, $f_S(P)$ is the 
number of chains in $P$ of the form $\hat{0}\prec x_1\prec x_2\prec\cdots\prec
x_k\prec\hat{1}$, where $\{\rank(x_i): 1\le i\le k\}=S$.
The $2^n$-tuple of flag numbers $f_S(P)$ (as $S$ ranges over all subsets of 
$[1,n]$) is called the {\em flag vector} of $P$.
The {\em flag $h$-vector} is obtained by performing inclusion-exclusion on the
flag vector.  Thus $h_S=\sum_{T\subseteq S} (-1)^{|S\setminus T|} f_T$ or,
equivalently, $f_S=\sum_{T\subseteq S} h_T$.
Write a generating function in noncommuting variables,
$\Psi(a,b)=\sum h_Su_S$, where $u_S=u_1u_2 \cdots u_n$ with $u_i=a$ if 
$i\not\in S$ and $u_i=b$ if $i\in S$.
For every Eulerian poset, there is a polynomial $\Phi(c,d)$ in 
noncommuting
variables $c$ and $d$ for which $\Psi(a,b) = \Phi(a+b,ab+ba)$.
The polynomial $\Phi(c,d)$ (or $\Phi_P(c,d)$ when we need to specify the 
poset $P$) is called the {\em $cd$-index} of the poset.
The $cd$-index of the dual of the Eulerian poset $P$ is obtained by reversing
every $cd$-word in the $cd$-index of $P$.
The coefficient of a $cd$-word $w$ is written as $[w]$ (or $[w]_P$).
We think of each $d$ as occupying two positions in a $cd$-word, namely,
the positions of $ab$ or $ba$ in the corresponding $ab$-words.
Let $\supp(w)$ be the set of positions of $d$ in $w$.

Stanley (\cite{Stanley-flag}) notes a useful variation of the $cd$-index.
The {\em $ce$-index} is obtained by replacing every $d$ in $\Phi(c,d)$ by
$(cc-ee)/2$.  Alternatively, one gets the $ce$-index from 
$\Psi(a,b)$---even for non-Eulerian posets---by
letting $c=a+b$ and $e=a-b$.
The $ce$-index is thus a polynomial in the noncommuting variables $c$ and $e$,
where for Eulerian posets the $e$'s occur only in pairs.
Write $L_Q$ for the coefficient of the word $v_Q=v_1v_2\cdots v_n$,
where $v_i=c$ if $i\not\in Q$ and $v_i=e$ if $i\in Q$.
The vector of coefficients of the $ce$-index, $(L_Q(P))$, is also known as
the {\em $L$-vector} of $P$.

For an Eulerian poset $P$, $L_Q(P)=0$ unless $Q$ is an {\em even}
set, that is, $Q$ is the union of disjoint intervals of even cardinality.
We say 
{\em $Q$ evenly contains~$S$}, written $S\subseteq_e Q$, if
$S$ and $Q$ are even sets, 
$S\subseteq Q$, 
and the difference set $Q\setminus S$ is 
also an even set.
An ``Eulerian'' $ce$-word $v_Q$ is converted to a sum of $cd$-words by
replacing consecutive pairs of $e$'s in $v_Q$ by $cc-2d$ so that no $e$'s
remain.
This means that a $cd$-word $w$ occurs in the expansion of a $ce$-word $v_Q$
if and only if $\supp(w)\subseteq_e Q$.
Thus the coefficient in the $cd$-index of a $cd$-word $w$ in which $d$ occurs 
$r$ times is 
\begin{equation}
\label{cdce}
[w]=(-2)^r\sum_{\supp(w)\subseteq_e Q} L_Q.
\end{equation}
(See \cite{Bayer-Hetyei1} for more information on $L$-vectors.)

In determining the cone of flag vectors of all graded posets 
(\cite{Billera-Hetyei}), Billera and Hetyei construct sequences of posets
with convergent (normalized) flag vectors.
Bayer and Hetyei apply a doubling operation to some of these to get 
sequences of Eulerian posets. 
Given an interval $I=[i,j]\subseteq [1,n]$, a rank $n+1$ poset $P$ and a
positive integer $N$, let $D^N_I(P)$ be the rank $n+1$ poset obtained by 
replacing
$P_I$, the subposet of $P$ consisting of elements with ranks in $I$,
by $N$ copies of itself.
The (horizontal) {\em double} $DP$ of a poset $P$ is the result of starting with
$P$ and successively applying the operators $D^2_{\{i\}}$, for $1\le i\le n$.
(In the Hasse diagram of $P$ every edge is replaced by $\Join$.)
For \iv\ a set of subintervals of $[1,n]$, \iv\ is an {\em even interval
system} if (1)~no interval of \iv\ is contained in another, (2)~every
interval of \iv\ is of even cardinality, and (3)~the intersection of any
two intervals of \iv\ is of even cardinality.
For each even interval system \iv\ over $[1,n]$, there exists a sequence of 
Eulerian posets, $DP(n,\iv, N)$, whose normalized flag vectors (and hence,
normalized $cd$-indices and $ce$-indices) converge.
These are obtained by starting with a rank $n+1$ chain, successively
applying the operators $D^N_I$ for the intervals $I\in \iv$, and finally
taking the horizontal double. 

For \iv\ an interval system of $k$ intervals, 
write $L_S(DP(n,\iv))={}$ \linebreak 
$\lim_{N\rightarrow\infty} L_S(DP(n,\iv,N))/N^k$.
(Here $2^n N^k$ is the number of maximal chains in $DP(n,\iv,N)$.)
The symbol $DP(n,\iv)$ is referred to as a {\em limit poset}.
These $ce$-index coefficients are given by the formula
\begin{equation}
\label{L}
 L_S(D P(n,{\mathcal I}))
=\sum _{j=0}^k (-1)^j                                 
\left|\left\{1\leq i_1<\cdots <i_j\leq k :\:
I_{i_1}\cup \cdots \cup I_{i_j}=S
\right\}\right|, 
\end{equation}
where $\iv=\{I_1,I_2,\ldots,I_k\}$.                          
See \cite{Bayer-Hetyei1} for details.
The formula applies for non-Eulerian limit posets as well; in that case
it can give nonzero $L_Q$ for noneven sets $Q$.

We use one other result stated and proved in \cite{Bayer-Hetyei1} (but implicit
in \cite{Billera-Liu}).
\begin{proposition}[Inequality Lemma]
\label{ineqlemma}
Let $T$ and $V$ be subsets of $[1,n]$ such that for every maximal interval $I$
of $V$, $|I\cap T|\le 1$.
Write $S=[1,n]\setminus V$.
For $P$ any rank $n+1$ Eulerian poset,
\[\sum_{R\subseteq T}(-2)^{|T\setminus R| } f_{S\cup R}(P)\ge 0. \]
Equivalently,
\[ (-1)^{|T|}\sum_{T\subseteq Q\subseteq V} L_Q(P)\ge 0.\] 
\end{proposition}

\section{The Main Result}

\begin{theorem} \mbox{ }
\label{mainth}
\begin{enumerate}
\item \label{th_pt1}
For the following $cd$-words $w$, the coefficient of $w$ as a function
of Eulerian posets has greatest lower bound 0 and has no upper bound:
\begin{enumerate}
\item $c^idc^j$, with $\min\{i,j\}\le 1$
\item $c^idcd\cdots cdc^j$ (at least two $d$'s alternating with $c$'s, $i$ and
      $j$ unrestricted)
\end{enumerate}
\item \label{th_pt2} The coefficient of $c^n$ in the $cd$-index of every 
      Eulerian poset is 1.
\item \label{th_pt3}
For all other $cd$-words $w$, the coefficient of $w$ as a function of
Eulerian posets has neither lower nor upper bound.
\end{enumerate}
\end{theorem}
Note.  For $n\ge 5$, there are $\lfloor {n-2\choose 2}/3 \rfloor + 4$ $cd$-words
of the types described in Part~\ref{th_pt1}.
This is a small portion of the $cd$ words for large $n$.

\vspace*{6pt}

\begin{proof}
The fact that the coefficient of $c^n$ is 1 is immediate from the definition
and is included only for completeness.

Let $w$ be any $cd$-word containing $r$ copies of $d$, with $r\ge 1$.
Let \iv\ be the set of two-element intervals of the positions of $d$ in $w$.
Compute the coefficient of $w$ in the $cd$-index of $DP(n,\iv)$.
If $Q$ properly contains $\supp(w)$, then by equation~(\ref{L}), 
$L_Q(DP(n,\iv))=0$.
So by equations~(\ref{cdce}) and~(\ref{L}), for $DP(n,\iv)$, the 
coefficient of $w$ is $[w]=(-2)^rL_{\supp(w)}(DP(N,\iv))=(-2)^r(-1)^r=2^r$.
This is the limit as $N$ goes to infinity of $1/N^r$ times $[w]$ for 
$DP(n,\iv,N)$.
So $(DP(n,\iv,N))$ is a sequence of Eulerian posets with $cd$-coefficients
$[w]$ not bounded above.

To show nonnegativity in Part~\ref{th_pt1}, we use equation~(\ref{cdce}) and
the Inequality Lemma.
If $w=dc^j$, let $S=\emptyset$ (so $V=[1,n]$) and $T=\{1\}$.
Then the coefficient of $w$ is $[w]= 2 (-1)\sum_{T\subseteq Q\subseteq V}
L_Q\ge 0$.
If $w=cdc^j$ let $S=\{1\}$ (so $V=[2,n]$) and  $T=\{2\}$.
Then the coefficient of $w$ is $[w]= (-2)\sum_{\{2,3\}\subseteq_e Q} L_Q
={}$ \linebreak $2(-1)\sum_{T\subseteq Q\subseteq V} L_Q\ge 0$, because $L_Q$ is
zero unless $Q$ is an even set.
The cases of $w=c^id$ and $w=c^idc$ follow by duality.

Let $w=c^idcd\cdots cdc^j$, with $d$ occurring $r$ times, $r\ge 2$; thus
$\supp(w)=\{i+1,i+2,i+4,i+5,\ldots,i+3r-2,i+3r-1\}$.
Let $S=\{i+3,i+6, \ldots, i+3r-3\}$, $V=[1,n]\setminus S$ and 
$T=\{i+2,i+4,i+7, \ldots, i+3r-5, i+3r-2\}$.
Here $S$ is the set of positions of the $c$'s between $d$'s and $T$
is a set of one position for each $d$, adjacent to the positions of the interior
$c$'s.
The coefficient of $w$ is $[w]= (-2)^r\sum_{\supp(w)\subseteq_e Q} L_Q$.
The set $Q$ evenly contains $\supp(w)$ if and only if $Q$ is an even set
and $T\subseteq Q\subseteq V$.
Since  $L_Q=0$ unless $Q$ is an even set,
$[w]=2^r(-1)^r\sum_{T\subseteq Q\subseteq V} L_Q \ge 0$.

The double of the chain, $DC^{n+1}$, has $cd$-index $c^n$, so for the 
$cd$-coefficients in Part~I, the lower bound of 0 is actually attained.

It remains to show that the coefficients of the $cd$ words in Part~\ref{th_pt3}
can be arbitrarily negative.  We use several lemmas.

\begin{lemma}
\label{even}
For every even $n\ge 4$ the coefficient of $dc^{n-4}d$ as a function of
Eulerian posets has no lower bound.
\end{lemma}

\begin{proof}
Let  $\iv=\{[1,n]\}$.
By equation~(\ref{L}) the only nonzero entries in the $L$-vector of 
$DP(n,\iv)$ are $L_\emptyset=1$ and $L_{[1,n]}=-1$.
By (\ref{cdce}) the coefficient of $dc^{n-4}d$ in the $cd$-index of 
$DP(n,\iv)$ is $(-2)^2(-1)=-4$.
This is the limit as $N$ goes to infinity of $1/N^2$ times $[dc^{n-4}d]$ for 
$DP(n,\iv,N)$.
So $(DP(n,\iv,N))$ is a sequence of Eulerian posets with $cd$-coefficients
$[dc^{n-4}d]$ not bounded below.
(A formula of Ehrenborg and Readdy (\cite{Ehrenborg-Readdy}) gives directly
that the $cd$-index of $DP(n,\iv,N)$ is $(N+1)c^n-N(cc-2d)^{n/2}$.)
\end{proof}

In \cite{Bayer-Hetyei2} Bayer and Hetyei discuss constructions of Eulerian
posets whose normalized $L$-vectors converge to sums of $L$-vectors of 
non-Eulerian Billera-Hetyei limit posets.
(A few examples are found in \cite[Appendix A]{Bayer-Hetyei1}.)

\begin{lemma}
\label{odd}
For every odd $n\ge 7$ the coefficient of $dc^{n-4}d$ as a function of 
Eulerian posets has no lower bound.
\end{lemma}

\begin{proof}
Write $C^{n+1}$ for the chain of rank $n+1$.
Let 
\[P^I(N)=
  D^{N+1}_{[1,2]}D^{N+1}_{[3,n-3]}D^{N+1}_{[4,n-2]}D^{N+1}_{[n-1,n]}(C^{n+1});\]
let 
\[P^{II}(N)=D^{N+1}_{[1,n-3]}D^{N^2}_{[3,n-2]}D^{N+1}_{[4,n]}(C^{n+1});\] and
let 
\[P^{III}(N)=D^{N^4}_{[1,n]}(C^{n+1}).\]
Create a poset $P(N)$ from these three posets by identifying
the elements of $P^{II}(N)$ with the elements of $P^I(N)$ at ranks 0, 
1, 2, $n-1$, $n$, and $n+1$, and then identifying 
the elements of $P^{III}(N)$ with the elements of $P^I(N)$ and 
$P^{II}(N)$ only at ranks 0 and $n+1$.
The doubles $DP(N)$ of these posets are Eulerian, and the normalized $L$-vectors
converge as $N$ goes to infinity.
Write $L_Q(DP)=\lim_{N\rightarrow \infty} L_Q(DP(N))/f_{[1,n]}(DP(N))$.
Then $L_Q(DP)=L_Q(DP(n,\iv_1)) + L_Q(DP(n,\iv_2)) + L_Q(DP(n,\iv_3))$,
where $\iv_1=\{[1,2],[3,n-3],[4,n-2],[n-1,n]\}$,
$\iv_2=\{[1,n-3],[3,n-2],[4,n]\}$, and $\iv_3=\{[1,n]\}$.
The only nonzero $L_Q$ for which $\{1,2,n-1,n\}\subseteq_e Q$
are $L_{\{1,2,n-1,n\}}(DP)=1$, $L_{[1,n]\setminus\{3\}}(DP)=-1$ and
$L_{[1,n]\setminus\{n-2\}}(DP)=-1$, so by equation~(\ref{cdce})
the coefficient of $dc^{n-4}d$ in the $cd$-index of $DP$ is $-4$.
This is the limit as $N$ goes to infinity of $1/f_{[1,n]}(DP(N))$ times 
$[dc^{n-4}d]$ for $DP(N)$.
So $(DP(N))$ is a sequence of Eulerian posets with $cd$-coefficients
$[dc^{n-4}d]$ not bounded below.
(In fact, a flag vector calculation gives the coefficient of $dc^{n-4}d$
for $DP(N)$ as $4(N^2-N^4)$.)
\end{proof}

The proof of Lemma~2 asserts that $DP(N)$ is Eulerian.
It is easy to check by equation~(\ref{L}) that 
$L_Q(DP(n,\iv_1)) + L_Q(DP(n,\iv_2)) + L_Q(DP(n,\iv_3))=0$ if $Q$ is
not an even set.
This condition must hold if every $DP(N)$ is an Eulerian poset.
But to prove that $DP(N)$ is Eulerian requires us to show that
every interval of the poset has the same number of elements of even rank and
of odd rank.
We show the details in one case.
Let $[x,y]$ be an interval of $P(N)$ with $x$ of rank~2 and $y$ of rank~$n-1$.
For the Eulerian condition to hold on corresponding intervals in $DP(N)$, the
interval $[x,y]$ of $P(N)$ must have one more element of even rank than of 
odd rank.
If $x$ and $y$ are in the subposet $P^{III}(N)$, then $[x,y]$ has exactly
one element of each rank, so the condition is met.
Suppose $x$ and $y$ are identified elements of $P^I(N)$ and $P^{II}(N)$.
In the open interval $(x,y)$ in $P^I(N)$, ranks~3 and~$n-2$ each have $N+1$
elements and each other rank has $(N+1)^2$ elements.
In the open interval $(x,y)$ in $P^{II}(N)$, each rank has $N^2$ elements.
So the number of even-rank elements in $[x,y]$ is $2 + ((N+1)^2+N^2)(n-5)/2$,
and the number of odd-rank elements in $[x,y]$ is 
$2(N+1)+2N^2+((N+1)^2+N^2)(n-7)/2$.
The difference is 1.
Note that neither $P^I(N)$ nor $P^{II}(N)$ satisfies the Eulerian condition
for $[x,y]$ by itself.
The two subposets balance each other to achieve the Eulerian property.
This works for all intervals.

\begin{lemma}
\label{ccdcc}
The coefficient of $ccdcc$ as a function of rank 7 Eulerian posets has no 
lower bound. 
\end{lemma}

\begin{proof}
The following limit poset is given in 
Appendix~A of \cite{Bayer-Hetyei1}.
Let $P^I(N)=D^N_{[1,2]}D^N_{[2,6]}(C^7)$ and 
$P^{II}(N)=D^N_{[1,5]}D^N_{[5,6]}(C^7)$.
Let $P(N)$ be formed from these two posets by identifying the elements at
ranks 0, 1, 6, and 7.
The double $DP(N)$ of this poset is Eulerian.
In the limit,
the normalized $L$-vector includes the following values:
$L_{34}= L_{1234}= L_{3456}=0$, and $L_{123456}=1$.
These are the $L_Q$ that contribute to the coefficient of $ccdcc$ in the 
$cd$-index, $[ccdcc]=
-2(L_{34} + L_{1234} + L_{3456} +  L_{123456})=-2$.
As argued before, 
this gives a sequence of Eulerian posets with $cd$-coefficients $[ccdcc]$ not
bounded below.
(In fact, for $DP(N)$, $[ccdcc]=-2(N-1)^2$.)
\end{proof}

\begin{lemma}
\label{extend}
Let $u$ and $v$ be $cd$-words.  
If the coefficient of $u$ as a function of Eulerian posets
has no lower bound, then the coefficients of $uv$ and $vu$ as functions of
Eulerian posets have no lower bounds.
\end{lemma}
\begin{proof}
In \cite{Stanley-flag} Stanley considers a ``join'' operation, which produces
an Eulerian poset $P*Q$ from two Eulerian posets $P$ and $Q$.
He shows that the $cd$-indices satisfy $\Phi_{P*Q}(c,d)=\Phi_P(c,d)\Phi_Q(c,d)$.
Let $u$ be a $cd$-word of length $m$ and $v$ a $cd$-word of length $n$.
Let $B$ be the rank $n+1$ Boolean algebra.
Every $cd$-word of length $n$ has a positive coefficient in the $cd$-index of
$B$.
(This is proved most easily from the Ehrenborg-Readdy
formula for the $cd$-index of a pyramid  in \cite{Ehrenborg-Readdy}.) 
Let $P_N$ be a sequence of rank $m+1$ Eulerian posets for which 
$\displaystyle\lim_{N\rightarrow \infty}[u]_{P_N}=-\infty$.
Then 
$\displaystyle\lim_{N\rightarrow \infty}[uv]_{P_N*B}= \lim_{N\rightarrow \infty}[vu]_{B*P_N}=
-\infty$.
\end{proof}

We now complete the proof of the theorem.
Every $cd$-word not included in Parts~1 and~2 of the theorem contains the
subword $ccdcc$ or a subword of the form $dc^{n-4}d$ for $n-4\ne 1$.
Thus, by Lemmas~1 through~4, the coefficients of these $cd$-words 
as functions of Eulerian posets have no lower bounds.
\end{proof}

\noindent {\bf Acknowledgments:}  
The referee was most generous with
suggestions for improving the paper.
I also wish to thank G\'{a}bor Hetyei for 
introducing me to the construction of Eulerian posets crucial to this paper,
and for other helpful discussions.

\providecommand{\bysame}{\leavevmode\hbox to3em{\hrulefill}\thinspace}

\end{document}